\documentclass{SCAEOL}
\numberwithin{equation}{section}
\usepackage{algorithm}
\usepackage{algorithmic}
\usepackage{xypic}
\usepackage{multirow}

\begin{document}

\Year{2017} %
\Month{January}
\Vol{60} %
\No{6} %
\BeginPage{1} %
\EndPage{XX} %
\AuthorMark{Huang D D {\it et al.}}
\ReceivedDay{March 14, 2016}
\AcceptedDay{March 22, 2017}
\DOI{10.1007/s11425-000-0000-0} 

\title[computing the factor ring of an ideal]{Algorithm for computing the factor ring of an ideal in Dedekind domain with finite rank}{}


\author[1]{Dandan Huang}{Corresponding author}
\author[2,3]{Yingpu Deng}{}

\address[{\rm1}]{Department of Cyber Space Security, School of Software Engineering,}
\address{Jinling Institute of Technology, Nanjing {\rm211169}, P.R. China;}

\address[{\rm2}]{Key Laboratory of Mathematics Mechanization, NCMIS, Academy of Mathematics and Systems Science,}
\address{Chinese Academy of Sciences, Beijing {\rm100190}, P.R. China}

\address[{\rm3}]{University of Chinese Academy of Sciences, Beijing {\rm100049}, P.R. China}

\Emails{huangdd@jit.edu.cn,
dengyp@amss.ac.cn}\maketitle


 {\begin{center}
\parbox{14.5cm}{\begin{abstract}
We give an algorithm for computing the factor ring of a given ideal
in a Dedekind domain with finite rank, which runs in deterministic and polynomial-time.
We provide two applications of the algorithm: judging whether a given ideal is prime or prime power.
The main algorithm is based on basis representation
of finite rings which is computed via Hermite and
Smith normal forms.\vspace{-3mm}
\end{abstract}}\end{center}}



\renewcommand{\baselinestretch}{1.2}
\begin{center} \renewcommand{\arraystretch}{1.5}
{\begin{tabular}{lp{0.8\textwidth}} \hline \scriptsize
{\bf Citation:}\!\!\!\!&\scriptsize Huang D D, Deng Y P. \makeatletter\@titlehead.
Sci China Math, 2017, 60,
 doi:~\@DOI\makeatother\vspace{1mm}
\\
\hline
\end{tabular}}\end{center}

\baselineskip 11pt\parindent=10.8pt  \wuhao
\section{Introduction}

Computing invariants of algebraic number fields, such as integral bases, discriminants and ideal class groups, is important both for its own sake and for its numerous applications. The practical completion of this task, which is usually called Dedekind program, has been one of the major achievements of computational number theory in recent years by many people, especially studied by Cohen in his books \cite{cohen,coh}.

In this paper, we address a fundamental problem of computing the
factor ring of a given ideal in Dedekind domain of finite rank,
where finite rank means that the ring as a $\mathbb{Z}-$module is finitely generated.
The main contribution of this paper is to give a deterministic polynomial-time algorithm, we call it
the main algorithm, which outputs a basis representation of the factor ring.
The concept of basis representation was first proposed
by Lenstra \cite{len} for describing finite fields.
Kayal et al. stated the formal definition of basis
representation for finite rings in \cite{KS}.
The idea of computing basis representation of factor rings
does not appear in the algorithms of \cite{cohen,coh} for dealing with Dedekind program.

We construct two important applications connected with the main algorithm in this paper.
In Dedekind domain, every nonzero ideal can
be written as a product of prime ideals in a unique way.
Like the famous problem of primality testing, one might to judge whether a given ideal of
Dedekind domain is prime or not. Actually,
the first application of this paper is to determine the primality of nonzero ideals in Dedekind domain with finite rank. The other one is to judge whether a given ideal is prime ideal power or not.

In \cite{cohen,coh}, Cohen stated an algorithm for
deciding the primality of nonzero ideals in the ring of algebraic
integers of a number filed, which is also Dedekind domain of finite rank.
Since the algorithm uses the factorization of univariate polynomials over finite fields,
which can not be computed in deterministic and polynomial complexity.
The algorithm of \cite{cohen,coh} does not run in
deterministic and polynomial-time.

By the fact that a nonzero proper ideal of Dedekind domain
is prime if and only if the corresponding factor ring is a field, we can apply the main algorithm
to deduce the first application for testing prime ideals.
Moreover, our prime ideal test runs in deterministic and polynomial-time.
The important fact we used in the analysis of computational complexity is
 the fact of field testing in \cite{ADM}, which stated that field testing of
finite rings in basis representation is a \textbf{P} problem.
We will analyze the asymptotical time bound of the field testing algorithm of \cite{ADM} in this paper.

Besides, due to the fact that a nonzero proper ideal of Dedekind domain
is a prime power if and only if its factor ring is a local ring,
we obtain the second application of the main algorithm.
The prime ideal power test of this application uses the local ring test
which is presented by Staromiejski in \cite{sta}.
The computational complexity of the prime ideal power test is also deterministic and polynomial.

The paper is organized as follows. In Section \ref{SecPro}, we
state the core problem precisely and recall some basic definitions, such as Dedekind domain and basis representation.
In Section \ref{SecMain}, we describe the main algorithm of this paper explicitly.
We prove the correctness of the main algorithm and analyse its computational complexity in Section \ref{Seccom}.
In Section \ref{Secapp}, we provide two immediate applications of the main algorithm.
And in Section \ref{Secimp}, we illustrate several different examples to show the implementation aspects of all the algorithms.
Section \ref{Seccon} is devoted to conclusions on the relevant analysis of our algorithms.

Throughout the paper, all rings are assumed to be commutative and
with multiplicative identity, written as $1$, and $1\neq0$. We
denote by $\mathrm{M}(t)$ an upper bound for the number of bit
operations required to multiply two $\lceil t\rceil$ bit integers.
By a result of \cite{sst},
$\mathrm{M}(t)=O(t\mathrm{log}t\mathrm{log}\mathrm{log}t)$.
Similarly, by $\mathrm{B}(t)$ we denote the number of bit operations
of the operation which is the application of the Chinese remainder
theorem with moduli consisting of all primes less than $t$. We can
take $\mathrm{B}(t)=O(\mathrm{M}(t)\mathrm{log}t)$ according to \cite{hmc}.
We denote by $\omega$ the exponent for matrix multiplication, and
$2<\omega\leq3$.

\section{Problem}\label{SecPro}

We begin with recalling the definition of Dedekind
domain, which can be found in the book \cite{Jan}.

\begin{definition}[Dedekind Domain]\quad
A ring $\mathcal{O}$ is a Dedekind domain if it is a noetherian
integral domain such that the localization
$\mathcal{O}_{\mathfrak{p}}$ is a discrete valuation ring for every
nonzero prime ideal $\mathfrak{p}$ of $\mathcal{O}$.
Moreover, if $\mathcal{O}$ as a $\mathbb{Z}-$module is of finite rank, then we call $\mathcal{O}$
a Dedekind domain with finite rank.
\end{definition}

Notice that every ideal in a Dedekind domain can be generated by at most two elements.
We will show some concrete examples of Dedekind domains in the following.

\begin{example}\quad
\begin{enumerate}
\item
The Gaussian Domain $\mathbb{Z}[i]$ is the ring of integers of the quartic cyclotomic field $\mathbb{Q}(i)$,  where $i=\sqrt{-1}$, which is a Dedekind domain of rank $2$.
\item
The Eisenstein Domain $\mathbb{Z}[\omega]$ is the ring of integers of the cubic cyclotomic field
$\mathbb{Q}(\omega)$, where $\omega=(-1+\sqrt{-3})/2$, which is also a Dedekind domain of rank $2$.
\item
The ring $\mathcal{O}=\mathbb{Z}[\sqrt[3]{2}]$ is a Dedekind domain of rank $3$.
\end{enumerate}
\end{example}

Actually, all the rings of algebraic integers of number fields are
Dedekind domains with finite rank. The ring of integers $\mathbb{Z}$ is a trivial Dedekind domain of rank $1$.

Throughout this paper, we focus on solving the following crucial problem to obtain the factor ring of a given ideal
in a Dedekind domain with finite rank.

\begin{problem}\label{mainproblem}\quad
In: $\mathcal{O}$, a Dedekind domain with finite rank; $I$, a nonzero ideal of $\mathcal{O}$.
Out: $R=\mathcal{O}/I$, the factor ring of $I$.
\end{problem}

To state the core problem (Problem \ref{mainproblem}) as an algorithmic problem, one needs to choose finite representations of
the input and the output. Usually, $\mathcal{O}$ is represented as a $\mathbb{Z}$-module, that is $\mathcal{O}=\mathbb{Z}\omega_1\oplus\ldots\oplus\mathbb{Z}\omega_n$, where
$\mathcal{W}=\{\omega_1,\ldots,\omega_n\}$ is called a $\mathbb{Z}$-basis of
$\mathcal{O}$. Besides, a multiplication table of $\mathcal{W}$ is given as
a sequence of integers $((c_{ijk})_{i,j,k=1,\ldots,n})$ such that
$$\omega_i\omega_j=\sum\limits_{k=1}^{n}c_{ijk}\omega_k.$$
Notice that, $c_{ijk}=c_{jik}$ because of $\omega_i\omega_j=\omega_j\omega_i$, for all $i, j, k$.

The output of Problem \ref{mainproblem} is a basis representation of the factor ring,
which is a finite ring under this case.
The definition of basis representation of finite rings is stated as follows by \cite{KS,len}.

\begin{definition}[Basis Representation]\label{basis}\quad
Let $R$ be a finite ring, a basis representation of $R$  is a
sequence of integers
$(m;d_1,\ldots,d_m;(l_{ijk})_{i,j,k=1,\ldots,m})$, where $m>0,\;
d_i\geq2$ and $0\leq l_{ijk}<d_k$, such that

$\mathrm{(1)}$ the additive group
$(R,+)=\mathbb{Z}_{d_1}\upsilon_1\oplus\ldots\oplus\mathbb{Z}_{d_m}\upsilon_m$,
where $d_i$ are the additive orders of generator $\upsilon_i$, and

$\mathrm{(2)}$ the  multiplication of
$\{\upsilon_1,\ldots,\upsilon_m\}$ is given by
$$\upsilon_i\upsilon_j=\sum\limits_{k=1}^{m}l_{ijk}\upsilon_k.$$

Integers $l_{ijk}$ are called structure constants.
\end{definition}

Now we restate Problem \ref{mainproblem} in the following way of representations of
the input and the output.

\smallskip
{\bf Representation of Problem \ref{mainproblem}.}\quad
In: $\mathcal{O}=\mathbb{Z}\omega_1\oplus\ldots\oplus\mathbb{Z}\omega_n$,
a multiplication table of $\mathcal{W}$ related to $\mathcal{O}$ is $((c_{ijk})_{i,j,k=1,\ldots,n})$;
$I=(\alpha,\beta)$, where $\alpha=\sum\limits_{i=1}^{n}a_i\omega_i$,
$\beta=\sum\limits_{i=1}^{n}b_i\omega_i$.\\
Out: $R=\mathcal{O}/I=(m;d_1,\ldots,d_m;(l_{ijk})_{i,j,k=1,\ldots,m})$.
\smallskip

We illustrate several explicit examples to explain the
above representations in detail.

\begin{example}\quad
In: $\mathcal{O}=\mathbb{Z}$, $I=(7)=7\mathbb{Z}$. \\
Out: $R=\mathcal{O}/I=(m;d_1,\ldots,d_m;(l_{ijk})_{i,j,k=1,\ldots,m})$.
\end{example}

Since this is a trivial example, we can compute a basis representation of the output immediately,
that is  $R=(1;7;(1))$.

\begin{example}\quad
In: $\mathcal{O}=\mathbb{Z}[\delta]$, $I=(5,2+\delta)$, where $\delta=\sqrt[3]{2}$.\\
Out:  $R=\mathcal{O}/I=(m;d_1,\ldots,d_m;(l_{ijk})_{i,j,k=1,\ldots,m})$.
\end{example}

\begin{example}\label{exmm}\quad
In: $\mathcal{O}=\mathbb{Z}[\theta]$, $I=(\theta-2)$, where $\theta=(1+\sqrt{-23})/2$.\\
Out: $R=\mathcal{O}/I=(m;d_1,\ldots,d_m;(l_{ijk})_{i,j,k=1,\ldots,m})$.
\end{example}

\begin{example}\quad
In: $\mathcal{O}=\mathbb{Z}[\gamma]$, $I=(23)$, where $\gamma^3=\gamma+1$.\\
 Out: $R=\mathcal{O}/I=(m;d_1,\ldots,d_m;(l_{ijk})_{i,j,k=1,\ldots,m})$.
\end{example}

One may not deduce basis representation of the factor rings for the above three examples easily,
until he or she applies the main algorithm described in the section that followed.

\section{The main algorithm}\label{SecMain}

In this section, we mainly describe an explicit algorithm to solve Problem \ref{mainproblem}.
Let us start with an elementary definition.

\begin{definition}[Norm]\quad
Let $\mathcal{O}$ be a Dedekind domain with finite rank,
$I$ be a nonzero ideal of $\mathcal{O}$.
We define the norm of $I$, $\mathcal{N}(I)$, to be the order of the factor ring
$\mathcal{O}/I$, i.e. $\mathcal{N}(I)=|\mathcal{O}/I|$.
Particularly, we denote $\mathcal{N}(\alpha)=\mathcal{N}(I)$ where
$I=(\alpha)$ is a principal ideal.
\end{definition}

According to the above definition, we are easy to deduce the multiplicative property of norm:
$\mathcal{N}(IJ)=\mathcal{N}(I)\mathcal{N}(J)$, where $I,\;J$ are
 nonzero ideals of $\mathcal{O}$. Computation of the norm of an ideal
is contained in the main algorithm being discussed later.

\begin{example}\quad
\begin{enumerate}

\item
If $\mathcal{O}=\mathbb{Z}$, $I=(7)$, then $\mathcal{N}(I)=|\mathbb{Z}/7\mathbb{Z}|=7$.

\item
If $\mathcal{O}=\mathbb{Z}[\gamma]$, $I=(23)$, where $\gamma^3=\gamma+1$. Since we have
$\mathcal{O}=\mathbb{Z}\oplus\mathbb{Z}\gamma\oplus\mathbb{Z}\gamma^2$, then
$\mathcal{N}(I)=|\mathbb{Z}[\gamma]/23\mathbb{Z}[\gamma]|=23^3=12167$.
\end{enumerate}
\end{example}

Let $I$ be a nonzero ideal
of a Dedekind domain $\mathcal{O}$ with finite rank.
We start with introducing an auxiliary algorithm to pre-compute a positive integer $h_I$ related to $I$.
If $I=(\alpha)$ is principal, then $h_I$ could be chosen as a multiple of $\mathcal{N}(\alpha)$.
Otherwise $I=(\alpha,\beta)$, both $\alpha$ and $\beta$ nonzero, $h_I$ could be chosen as a
 multiple of $\mathcal{N}(\alpha)$ and $\mathcal{N}(\beta)$. We call $h_I$ a multiple-norm of $I$
in this paper. Note that $h_I$ is a multiple of $\mathcal{N}(I)$, since $(\alpha)\subseteq I$,
that is to say $\mathcal{N}(I)|\mathcal{N}(\alpha)$.

\begin{algorithm}[htb]
\caption{Sub-algorithm}\label{sub}
\begin{algorithmic}[1]
\REQUIRE \quad $\mathcal{O}=\mathbb{Z}\omega_1\oplus\ldots\oplus\mathbb{Z}\omega_n$, a multiplication table of $\mathcal{W}$ related to $\mathcal{O}$ is $((c_{ijk})_{i,j,k=1,\ldots,n})$;
$I=(\alpha,\beta)$, where $0\neq\alpha=\sum\limits_{i=1}^{n}a_i\omega_i$,
$\beta=\sum\limits_{i=1}^{n}b_i\omega_i$.
\ENSURE \quad $h_I$, a multiple-norm of $I$.

\STATE
Compute two $n\times n$ integral matrices $A=(a_{ij}),
B=(b_{ij})$, where
$$a_{ij}=\sum\limits_{k=1}^{n}a_kc_{kij},\; b_{ij}=\sum\limits_{k=1}^{n}b_kc_{kij},\; \mathrm {for\; all}\;1\leq i, j\leq n.$$

\STATE
If $\beta=0$, $h_I\leftarrow|\mathrm{det}(A)|$; otherwise $h_I\leftarrow|\mathrm{det}(AB)|$.
\end{algorithmic}
\end{algorithm}

Next we illustrate the sub-algorithm (Algorithm \ref{sub}) by implementing the following example.

\begin{example}\label{exm}\quad
In: $\mathcal{O}=\mathbb{Z}[\delta]$, $I=(5,2+\delta)$, where $\delta=\sqrt[3]{2}$.
Out: $h_I$.

First we obtain $\mathcal{O}=\mathbb{Z}\oplus\mathbb{Z}\delta\oplus\mathbb{Z}\delta^2$,
and the multiplication table of $\mathcal{W}=\{1,\delta,\delta^2\}$ is
$(c_{11k})_{k=1,2,3}=(1,0,0)$,
$(c_{12k})_{k=1,2,3}=(c_{21k})_{k=1,2,3}=(0,1,0)$,
$(c_{22k})_{k=1,2,3}=(c_{13k})_{k=1,2,3}=(c_{31k})_{k=1,2,3}=(0,0,1)$,
$(c_{23k})_{k=1,2,3}=(c_{32k})_{k=1,2,3}=(2,0,0)$,
$(c_{33k})_{k=1,2,3}=(0,2,0)$.

{\bf Step 1}\quad After some computations, the $3\times3$ integral matrices $A,\;B$ related to $I$ are

$$A=
\left(
\begin{array}{ccc}
5 & 0  & 0  \\
0  & 5  & 0 \\
0  & 0 &  5
\end{array}\right);\quad
B=
\left(
\begin{array}{ccc}
2 & 1  &  0 \\
0  &  2 & 1 \\
 2 & 0 & 2
\end{array}\right)$$

{\bf Step 2}\quad By computation the sub-algorithm outputs $h_I=|\mathrm{det}(AB)|=1250$.
\end{example}

The main algorithm for computing the factor ring of a given ideal in Dedekind domain with finite rank is stated as follows.

\begin{algorithm}[htb]
\caption{Main algorithm}\label{main}
\begin{algorithmic}[1]
\REQUIRE \quad  $\mathcal{O}=\mathbb{Z}\omega_1\oplus\ldots\oplus\mathbb{Z}\omega_n$,
a multiplication table of $\mathcal{W}$ related to $\mathcal{O}$ is $((c_{ijk})_{i,j,k=1,\ldots,n})$;
$I=(\alpha,\beta)$ with $0\neq\alpha=\sum\limits_{i=1}^{n}a_i\omega_i$,
$\beta=\sum\limits_{i=1}^{n}b_i\omega_i$; $h_I$, a multiple-norm of $I$.
\ENSURE \quad either $I=\mathcal{O}$ or $R=(m;d_1,\ldots,d_m;(l_{ijk})_{i,j,k=1,\ldots,m})$,
a basis representation  of $R=\mathcal{O}/I$.

\STATE
Compute the integral matrices $A=(a_{ij}),\;B=(b_{ij})$ such that
$$a_{ij}=\sum\limits_{k=1}^{n}a_kc_{kij},\; b_{ij}=\sum\limits_{k=1}^{n}b_kc_{kij},\; \mathrm {for\; all}\;1\leq i, j\leq n.$$

\STATE
$H_A\leftarrow$ the Hermite normal form of $A^T$, $H_B\leftarrow$ the Hermite normal form of $B^T$.

\STATE
If $\beta=0$, $H_M\leftarrow H_A$;
otherwise $\tilde{H}_M=(0\;H_M)\leftarrow$ the Hermite normal form of $M=(H_A,H_B)$.

\STATE
$S\leftarrow$ the Smith normal form of $H_M$, where $S=VH_MU$, $S=diag(d_1,\ldots,d_n)$,
and $d_{i+1}|d_i$ for all $1\leq i<n$.

\STATE
$\mathrm{det}(S)\leftarrow\prod\limits_{i=1}^{n}d_i$,
if $\mathrm{det}(S)=1$, then output $I=\mathcal{O}$ and STOP.

\STATE
$\tilde{V}\leftarrow$ the inverse matrix of $V$ over the ring $\mathbb{Z}_{h_I}$, compute the $n\times n$ integral matrices $\tilde{A_k}$ ($1\leq k\leq n$) such that
\begin{equation}\label{Ak}
\tilde{A_k}=(\tilde{a_{ijk}})_{1\leq i, j\leq
n}=(\tilde{V})^\mathrm{T}(c_{ijk})_{1\leq i, j\leq n}\tilde{V}.
\end{equation}

\STATE
For all $1\leq i, j\leq n$, compute
\begin{equation}\label{tij}
\left(
\begin{array}{c}
\tilde{t_{ij1}} \\
\vdots \\
\tilde{t_{ijn}}
\end{array}\right)=V
\left(
\begin{array}{c}
\tilde{a_{ij1}} \\
\vdots \\
\tilde{a_{ijn}}
\end{array}\right)
\end{equation}

\STATE
For $1\leq k\leq n$, $\pi_k\leftarrow$ the natural ring
homomorphism from $\mathbb{Z}$ to $\mathbb{Z}_{d_k}$ by
$\pi_k(a)=a\;\mathrm{mod}\;d_k$, $l_{ijk}\leftarrow\pi_k(\tilde{t_{ijk}})$,
then output $R=(m;d_1,\ldots,d_m;(l_{ijk})_{i,j,k=1,\ldots,m})$
such that all $d_1,\ldots,d_m$ are greater than $1$, where $1\leq m\leq n$.
\end{algorithmic}
\end{algorithm}

Note that $a\;\mathrm{mod}\;x$ denotes the smallest nonnegative residue of $a$
modulo $x$. As for the definitions of Hermite and Smith normal forms in
this paper, which originated in \cite{her} and \cite{smi} respectively, the reader can refer to the book \cite{cohen}.

\begin{remark}\quad
Under the pre-computation of $h_I$, a multiple-norm of $I$, we find that the main
algorithm (Algorithm \ref{main}) would be more practical, especially after we receive a
smaller $h_I$. Maybe we know such $h_I$ in advance, then the sub-algorithm
(Algorithm \ref{sub}) could be omitted.
\end{remark}

\section{Correctness and computational complexity}\label{Seccom}

The correctness and the computational complexity of the sub-algorithm and the main
algorithm are discussed in this section.

\begin{lemma}\label{ab}\quad
The sub-algorithm (Algorithm \rm{\ref{sub})} is correct.
\end{lemma}

\begin{proof}\quad The correctness of Algorithm \ref{sub} follows easily from
the fact $|\mathrm{det}(A)|=\mathcal{N}(\alpha)$ and $|\mathrm{det}(AB)|=\mathcal{N}(\alpha\beta)$ respectively.
Indeed, for all $1\leq i\leq n$ we have
\[
\begin{array}{ll}
\alpha\omega_i&=\sum\limits_{j=1}^{n}a_j\omega_i\omega_j
=\sum\limits_{j=1}^{n}a_j\sum\limits_{k=1}^{n}c_{ijk}\omega_k\\\\
&=\sum\limits_{k=1}^{n}(\sum\limits_{j=1}^{n}a_jc_{ijk})\omega_k
=\sum\limits_{j=1}^{n}(\sum\limits_{k=1}^{n}a_kc_{ikj})\omega_j\\\\
&=\sum\limits_{j=1}^{n}(\sum\limits_{k=1}^{n}a_kc_{kij})\omega_j.
\end{array}
\]

Similarly,
$$\beta\omega_i=\sum\limits_{j=1}^{n}b_j\omega_i\omega_j
=\sum\limits_{j=1}^{n}(\sum\limits_{k=1}^{n}b_kc_{kij})\omega_j.$$
That is to say
$$(\alpha\omega_1,\ldots,\alpha\omega_n)=(\omega_1,\ldots,\omega_n)A^\mathrm{T},\;
(\beta\omega_1,\ldots,\beta\omega_n)=(\omega_1,\ldots,\omega_n)B^\mathrm{T}.$$
Since the ideals
$$(\alpha)=\mathbb{Z}\alpha\omega_1\oplus\ldots\oplus\mathbb{Z}\alpha\omega_n,\;
(\beta)=\mathbb{Z}\beta\omega_1\oplus\ldots\oplus\mathbb{Z}\beta\omega_n,$$
and using the fact of \cite[Chapter 2,
Theorem 2.4.13]{cohen}, one can deduce that $\mathcal{N}(\alpha)=|\mathrm{det}(A)|$ and
$\mathcal{N}(\beta)=|\mathrm{det}(B)|$. This completes the proof.
\end{proof}

\begin{theorem}\label{core}\quad
The main algorithm (Algorithm \rm{\ref{main}}) is correct.
\end{theorem}

\begin{proof}\quad First we let
$$(\alpha_1,\ldots,\alpha_n)=(\omega_1,\ldots,\omega_n)H_A,\;
(\beta_1,\ldots,\beta_n)=(\omega_1,\ldots,\omega_n)H_B.$$
It follows from the proof of Lemma \ref{ab} that
$$(\alpha)=\mathbb{Z}\alpha_1\oplus\ldots\oplus\mathbb{Z}\alpha_n,\;
\mathrm{and} \;\;(\beta)=\mathbb{Z}\beta_1\oplus\ldots\oplus\mathbb{Z}\beta_n.$$

Since the columns of $M=(H_A,H_B)$, treated as the ordinates representation
with respect to $\mathcal{W}=\{\omega_1,\ldots,\omega_n\}$, generate the ideal
$I=\mathcal{O}\alpha+\mathcal{O}\beta$.
Let $(\gamma_1,\ldots,\gamma_n)=(\omega_1,\ldots,\omega_n)H_M$,
then one may obtain
$I=\mathbb{Z}\gamma_1\oplus\ldots\oplus\mathbb{Z}\gamma_n$.
Since the Smith normal form of $H_M$ is $S$,
we have $\mathrm{det}(S)=\mathrm{det}(H_M)=\mathcal{N}(I)$.
Clearly, if $\mathcal{N}(I)=\mathrm{det}(S)=1$, then
we get the trivial case $I=\mathcal{O}$.

If $\mathrm{det}(S)>1$, let $(\eta_1,\ldots,\eta_n)=(\omega_1,\ldots,\omega_n)V^{-1}$.
After some computations we have
$$(\gamma_1,\ldots,\gamma_n)U=(\eta_1,\ldots,\eta_n)S=(d_1\eta_1,\ldots,d_n\eta_n).$$
Since both the transforming matrices $U,\;V$ are unimodular matrices,
we may get
$$\mathcal{O}=\mathbb{Z}\eta_1\oplus\ldots\oplus\mathbb{Z}\eta_n,\;
\mathrm{and} \;\;
I=\mathbb{Z}(d_1\eta_1)\oplus\ldots\oplus\mathbb{Z}(d_n\eta_n).$$
It yields that
$(R,+)=(\mathcal{O}/I,+)=\mathbb{Z}_{d_1}\bar{\eta_1}\oplus\ldots\oplus\mathbb{Z}_{d_n}\bar{\eta_n}$,
where $\bar{\eta_i}$ denotes the coset $\eta_i+I$ belonged to the
factor ring $R$, $i=1,\ldots,n$.

Since
$(\eta_1,\ldots,\eta_n)=(\omega_1,\ldots,\omega_n)V^{-1}$ and under some
concrete computations we obtain the following equality of $n\times n$ integral matrices:

\begin{equation}\label{nij}\quad
(\eta_i\eta_j)_{1\leq i, j\leq n}=(V^{-1})^\mathrm{T}(\omega_i\omega_j)_{1\leq i, j\leq n}V^{-1}.
\end{equation}

Denote $n\times n$ integral matrices

\begin{equation}\label{Ak2}\quad
A_k=(a_{ijk})_{1\leq i, j\leq
n}=(V^{-1})^\mathrm{T}(c_{ijk})_{1\leq i, j\leq n}V^{-1},\;1\leq
k\leq n
\end{equation}
and compute the following vectors
\begin{equation}\label{tij2}\quad
\left(
\begin{array}{c}
t_{ij1} \\
\vdots \\
t_{ijn}
\end{array}\right)=V
\left(
\begin{array}{c}
a_{ij1} \\
\vdots \\
a_{ijn}
\end{array}\right)
\end{equation}
for all $1\leq i, j\leq n$. Then one can verify that
$\eta_i\eta_j=\sum\limits_{k=1}^{n}t_{ijk}\eta_k$ by the expressions (\ref{nij}), (\ref{Ak2}) and (\ref{tij2}) for all $1\leq i,j\leq n$.
Moreover, in the factor ring $R$ we have
$$\bar{\eta_i}\bar{\eta_j}=\overline{\eta_i\eta_j}=\overline{\sum\limits_{k=1}^{n}t_{ijk}\eta_k}
=\sum\limits_{k=1}^{n}\pi_k(t_{ijk})\bar{\eta_k}.$$

Hence the required structure constants of the basis representation of $R$ are contained in
these $\pi_k(t_{ijk})$. It suffices to show that $\pi_k(t_{ijk})=\pi_k(\tilde{t_{ijk}})$ for all $1\leq i,j,k\leq n$.
Indeed, let $\pi$ be the natural ring
homomorphism from $\mathbb{Z}$ to $\mathbb{Z}_h$ by
$\pi(a)=a\;\mathrm{mod}\;h$.  Since
$\prod\limits_{k=1}^{n}d_k=\mathrm{det}(S)=\mathcal{N}(I)|h$, there is a commutative
diagram linked $\pi$ to $\pi_k$:
$$\xymatrix{
  \mathbb{Z} \ar[d]_{\pi} \ar[r]^{\pi_k} &  \mathbb{Z}_{d_k}       \\
  \mathbb{Z}_h \ar[ur]_{\phi_k}                     }$$
where $\phi_k$ is the natural ring homomorphism from $\mathbb{Z}_h$
to $\mathbb{Z}_{d_k}$ by
$\phi_k(a\;\mathrm{mod}\;h)=a\;\mathrm{mod}\;d_k$. Comparing the expressions (\ref{Ak}) and (\ref{tij})
of Algorithm \ref{main} with
the expressions (\ref{Ak2}) and (\ref{tij2}), we note that the matrix $\tilde{V}$ need not be $V^{-1}$,
but we have $\pi(V^{-1})=\tilde{V}$, where $\pi$ acts on a matrix by
mapping on each entry of it. Finally we have
 $$\pi_k(t_{ijk})=\phi_k\circ\pi(t_{ijk})
=\phi_k\circ\pi(\tilde{t_{ijk}})=\pi_k(\tilde{t_{ijk}}),$$
for all $1\leq i,j,k\leq n$. This leads to the correctness of the main algorithm (Algorithm \ref{main}).
\end{proof}

Next we analyse the computational complexity of Algorithms
\ref{sub} and \ref{main}. Before that we introduce the
algorithms for computing the Hermite and Smith normal forms
of integral matrices. We require some notation. If $A$ is a matrix
 over $\mathbb{Z}$, we let $\mathcal{L}(A)$ denote the lattice generated by the columns
 of $A$, and let $\mathrm{det(}\mathcal{L}(A))$ denote the determinant of this lattice. What we need are the following two
algorithms that originated in the results of
\cite{hmc}. And the method of \cite{ili} is also used to
obtain Proposition \ref{SNFalg}.

\begin{proposition}\label{HNFalg}\quad
There exists a deterministic algorithm that receives as input an $n\times m$
integral matrix $A$ of rank $n$ and a positive integer $h$ that is a
multiple of $\mathrm{det(}\mathcal{L}(A))$, and produces as output
the Hermite normal form $H$ of $A$ such that $AU=H$, where $U$ is an
$m\times m$ unimodular matrix. The running time of the algorithm is
$O(mn\mathrm{B(log}T)+mn^2\mathrm{B(log}h))$ bit operations, if the
entries of $A$ are bounded in absolute value by $T$.
\end{proposition}

\begin{proposition}\label{SNFalg}\quad
There exists a deterministic algorithm that receives as input an $n\times n$
nonsingular integral matrix $B$ and a positive integer $h$ that is a
multiple of $\mathrm{det}(B)$, and produces as output the
Smith normal form $S$ of $B$ and the transforming matrices $U,\;V$ such
that $VBU=S$, where $U,\;V$ are $n\times n$ unimodular matrices. The
running time of the algorithm is
$O(n^2\mathrm{B(log}T)+n^3\mathrm{B(log}h)\mathrm{log}h)$ bit
operations, if the entries of $B$ are bounded in absolute value by $T$.
\end{proposition}

Now we compute the complexity of the sub-algorithm for outputting a multiple-norm of $I$.

\begin{lemma}\label{ab1}\quad
The time complexity of the sub-algorithm (Algorithm \rm{\ref{sub}})
is $O(n^3\mathrm{B}(n\mathrm{log}nT))$ bit operations, if all
integers $c_{ijk},\;a_i,\;b_i,\;1\leq i,j,k\leq n$ are bounded in
absolute value by $T$.
\end{lemma}

\begin{proof}\quad First, computing the matrices $A$ and $B$ of Step 1
in Algorithm \ref{sub} can be done in $O(n^3\mathrm{M(log}nT))$ bit operations.

Since $|a_{ij}|,\;|b_{ij}|\leq nT^2$,
 one may obtain that the entries of the matrix $AB$ are
bounded in absolute value by $n^3T^4$. By the Hadamard inequality,
we have $h_I\leq n^{n/2}(n^3T^4)^n=n^{7n/2}T^{4n}$.
Thus we may use small primes
modular computation to compute the determinants of $A$ and $B$ in Step 2.
That is to say, we first apply Gaussian elimination to compute the determinants of $A$
and $B$ modulo small primes $p$ no more than $t=O(n\mathrm{log}nT)$,
then recover $|\mathrm{det}(A)|$ and $|\mathrm{det}(B)|$ by the
Chinese remainder theorem (see \cite{hmc} for details). Hence it costs
$O(n^3\mathrm{B}(n\mathrm{log}nT))$ bit operations to obtain the
value of $h_I$. And the total complexity of Algorithm \ref{sub} is
$O(n^3\mathrm{B}(n\mathrm{log}nT))$ bit operations.
\end{proof}

The next is the analysis of the computational complexity of the main algorithm (Algorithm \ref{main}) for computing basis representation of the factor ring.

\begin{theorem}\label{core1}\quad
The time complexity of Algorithm \rm{\ref{main}} is $O(n^3\mathrm{B}(\mathrm{log}nT)+n^4\mathrm{B}(\mathrm{log}h_I)\mathrm{log}h_I)$
bit operations, where $T$ is as Lemma $\ref{ab1}$.
\end{theorem}

\begin{proof}\quad The time complexity of Step 1 in Algorithm \ref{main} is the same as the one in Lemma \ref{ab1}. In Step 2, it takes
$O(n^2\mathrm{B}(\mathrm{log}nT)+n^3\mathrm{B}(\mathrm{log}h_I))$ bit
operations to obtain the Hermite normal forms of
matrices $A^\mathrm{T}$ and $B^\mathrm{T}$ by applying Proposition $\ref{HNFalg}$.
Similarly, we apply Proposition $\ref{HNFalg}$ to $M$ in Step 3,
and obtain the Hermite normal form $\tilde{H}_M$
of $M$ in $O(n^3\mathrm{B}(\mathrm{log}h_I))$ bit operations, since the
entries of $M=(H_A,H_B)$ are bounded in absolute value by $h_I$.

By applying Proposition \ref{SNFalg} to $H_M$, we find that the time complexity
of computing the Smith normal form $S$ and the transforming matrices $U,\;V$
of $H_M$ is
$O(n^3\mathrm{B}(\mathrm{log}h_I)\mathrm{log}h_I)$ bit operations in Step 4 of Algorithm \ref{main}.

In the last three steps of Algorithm \ref{main}, it suffices to
compute all values of $\pi_k(\tilde{t_{ijk}})$ as follows.
Since $\mathrm{det}(V)=\pm1$, one can perform row reductions on $V$ to
compute the inverse matrix $\tilde{V}$ of $V$ over the ring
$\mathbb{Z}_{h_I}$, which can be done in $O(n^3\mathrm{log}^2h_I)$ bit
operations. Also one can obtain all values of
$\pi(\tilde{t_{ijk}}),\;1\leq i,j,k\leq n$, by computing the product
of matrices over $\mathbb{Z}_{h_I}$ in the expressions (\ref{Ak}) and (\ref{tij}) of Algorithm \ref{main}.
Then we calculate all
$\pi_k(\tilde{t_{ijk}})=\phi_k\circ\pi(\tilde{t_{ijk}})$, where $\phi_k$ and $\pi$ are well-defined in the proof of Theorem \ref{core}. All these computations can be done in
$O(n^4\mathrm{M}(\mathrm{log}h_I)+n^3\mathrm{log}^2h_I)=O(n^4\mathrm{log}^2h_I)$
bit operations.

Hence the total time complexity of Algorithm \ref{main} is
$O(n^3\mathrm{B}(\mathrm{log}nT)+n^4\mathrm{B}(\mathrm{log}h_I)\mathrm{log}h_I)$
bit operations, which is as asserted.
\end{proof}

\begin{remark}\label{size}\quad
If we take $h_I=\mathcal{N}(\alpha\beta)\leq
n^{7n/2}T^{4n}$ according to the sub-algorithm (Algorithm \ref{sub}), then
$\mathrm{log}h_I=O(n\mathrm{log}nT)$. So the input size of Algorithm \ref{main}
is $O(n^3\mathrm{log}T)$ bits. And the main algorithm is
deterministic and polynomial-time in the
input size for computing the factor ring of a given ideal in Dedekind domain with finite rank.
\end{remark}

\section{Applications of the main algorithm}\label{Secapp}

\subsection{Deciding whether a given ideal is prime}

First we apply the main algorithm to decide whether
a given nonzero ideal is prime in Dedekind domain with finite rank.
We begin with recalling the
algorithm of \cite{ADM} which states that field testing of
finite rings is of deterministic and
polynomial-time complexity. Now we exhibit
this algorithm in detail and show an explicit analysis of its computational complexity, which is
not analysed in \cite{ADM}.

\begin{algorithm}[htb]
\caption{Is-Field}\label{isfield}
\begin{algorithmic}[1]
\REQUIRE \quad $R=(m;d_1,\ldots,d_m;(l_{ijk})_{i,j,k=1,\ldots,m})$, a finite ring.
\ENSURE \quad TRUE iff $R$ is a field.

\STATE \label{isp}
If $d_1=\ldots=d_m$ is prime does not hold,
 return FALSE and STOP.

\STATE
If $m=1$, return TRUE and STOP.

\STATE\label{v1}
$p\leftarrow d_1$,  $f_1\leftarrow$ the minimal polynomial of the first generator
$\upsilon_1$ over $\mathbb{F}_{p}$.

\STATE \label{fre}
If $f_1$ is reducible over
$\mathbb{F}_{p}$, return FALSE and STOP.
\STATE
$m_1\leftarrow$ the degree of $f_1$,
if $m_1=m$,  return TRUE and STOP.

\FOR{$i=2$ to $m$}
\STATE\label{ui}
 $f_i\leftarrow$ the minimal polynomial of the $i-$th generator $\upsilon_i$ over
 $\mathbb{F}_{p}(\upsilon_1,\ldots,\upsilon_{i-1})$;
 \STATE
  $m_i\leftarrow$ the degree of $f_i$,
 if $f_i$ is reducible over $\mathbb{F}_{p}(\upsilon_1,\ldots,\upsilon_{i-1})$,\label{uii}
 return FALSE and STOP.
 \STATE
 If $\prod\limits_{j=1}^{i}m_j=m$, return TRUE and STOP.
\ENDFOR
\end{algorithmic}
\end{algorithm}

\begin{lemma}\label{field}\quad
Algorithm \rm{\ref{isfield}} is correct and runs in
$O(\mathrm{M(log}^{15/2}p)+m^6\mathrm{log}^3p)$ bit operations,
where $p=\mathrm{min}\{d_1,\ldots,d_m\}$.
\end{lemma}

\begin{proof}\quad The correctness follows from \cite{ADM}. We proceed with the proof of
the running time. The time complexity of Step \ref{isp} of Algorithm \ref{isfield} is dominated
by any known bound for deterministic primality testing. It
takes $O(\mathrm{M(log}^{15/2}p))$ bit operations by applying the AKS test of
\cite{AKS}. In Step \ref{v1},  computing the minimal polynomial of $\upsilon_1$  can be done in
$O(m^\omega\mathrm{log}m\mathrm{log}^2p)$ bit operations. Indeed, applying the method of \cite{sta} to
$R$, where
$(R,+)=\mathbb{F}_{p}\upsilon_1\oplus\ldots\oplus\mathbb{F}_{p}\upsilon_m$
is a $\mathbb{F}_{p}-$algebra, we can obtain the above complexity. It costs
$O(m^{(\omega+1)/2}\mathrm{log}m\mathrm{log}\mathrm{log}m\mathrm{log}^3p)$
bit operations to determine whether $f_1$ is reducible in Step \ref{fre}
(see \cite{sh} for details).

As to Step \ref{ui}, we describe a method for computing the minimal polynomial $f_i$
of $\upsilon_i$ over the field
$F_{i-1}=\mathbb{F}_{p}(\upsilon_1,\ldots,\upsilon_{i-1})$ firstly, where
$i>1$. It takes $O(m^\omega\mathrm{log}m\mathrm{log}^2p)$ bit
operations to compute a matrix $E\in\mathrm{M}_{(m+1)\times
m}(\mathbb{F}_{p})$ such that
\begin{equation}
\label{e}(1,\upsilon_i,\ldots,\upsilon_i^{m})=(\upsilon_1,\ldots,\upsilon_m)E^\mathrm{T}.
\end{equation}
The same technique as \cite{sta} is used for computing $E$. Then we calculate a great
linearly independent subset $S$ of
$\{\upsilon_1,\ldots,\upsilon_m\}$ over $F_{i-1}$ in the following
way. For instance, computing a great linearly independent subset of
$\{\upsilon_1,\upsilon_2\}$ is equivalent to solving the equation
(\ref{i-1}) of variables $x$ and $y$ belonged to $F_{i-1}$:
\begin{equation}\label{i-1}\quad
 x\upsilon_1+y\upsilon_2=0
\end{equation}
Since
$\{\upsilon_1^{t_1}\cdot\ldots\cdot\upsilon_{i-1}^{t_{i-1}}\;|\;0\leq
t_j<m_j,\; \mathrm{for}\; 1\leq j<i\}$ is a $\mathbb{F}_{p}-$basis
of $F_{i-1}$, we may write $x$ and $y$ in the coordinate
representation related to this $\mathbb{F}_{p}-$basis. Then (\ref{i-1}) is
converted into a linear system of equations over $\mathbb{F}_{p}$
with the help of the known structures of $F_{i-1}$ and $R$. The
system of equations owns $m$ equations and
$2\prod\limits_{j=1}^{i-1}m_j$ variables, where
$\prod\limits_{j=1}^{i-1}m_j\leq m$. It can be solved by performing
$2m^\omega$ operations in $\mathbb{F}_{p}$. We repeat this procedure
for $m-1$ steps by adding all generators
$\{\upsilon_2,\ldots,\upsilon_m\}$ one at a time to $\upsilon_1$,
then $S$ could be computed in
$(2+\ldots+m)m^\omega\mathrm{log}^2p=O(m^{\omega+2}\mathrm{log}^2p)$
bit operations. We may assume $S=\{\mu_1,\ldots,\mu_s\}$. As a
by-product one simultaneously receives a matrix $H\in\mathrm{
M}_{m\times s}(F_{i-1})$ such that

\begin{equation}\label{H}\quad
\left(
\begin{array}{c}
\upsilon_1 \\
\upsilon_2\\
\vdots\\
\upsilon_m
\end{array}\right)=H
\left(
\begin{array}{c}
\mu_1 \\
\vdots \\
\mu_s
\end{array}\right)
\end{equation}

Combining (\ref{e}) with (\ref{H}) we get
\begin{equation}\label{eh}\quad
\left(
\begin{array}{c}
_1 \\
\upsilon_i\\
\vdots\\
\upsilon_i^m
\end{array}\right)=EH
\left(
\begin{array}{c}
\mu_1 \\
\vdots \\
\mu_s
\end{array}\right)
\end{equation}

Finally one can perform row reductions on the matrix $EH$ to obtain the
minimal polynomial $f_i$ of $\upsilon_i$ over $F_{i-1}$. Hence
computing $f_i$ can be done in $O(m^{\omega+2}\mathrm{log}^2p)$ bit
operations. On the other hand, it costs
$O(m_i^{(\omega+1)/2}m^3\mathrm{log}m\mathrm{loglog}m\mathrm{log}^3p)$
bit operations for deciding whether $f_i$ is reducible or not in
Step \ref{uii} by using the method of \cite{sh},
where $m_i\leq m$. All in all, the total complexity of Algorithm
\ref{isfield} is $O(\mathrm{M(log}^{15/2}p)+m^6\mathrm{log}^3p)$ bit
operations, where $p=\mathrm{min}\{d_1,\ldots,d_m\}$.
\end{proof}

It is easy to deduce the algorithm for testing prime ideals in Dedekind domain with finite rank
by the main algorithm (Algorithm \ref{main}) and the field testing (Algorithm \ref{isfield}).
The prime ideal test is described as follows.

\begin{algorithm}[htb]
\caption{Is-Prime-Ideal}\label{primeideal}
\begin{algorithmic}[1]
\REQUIRE \quad $\mathcal{O}=\mathbb{Z}\omega_1\oplus\ldots\oplus\mathbb{Z}\omega_n$,
a multiplication table of $\mathcal{W}$ related to $\mathcal{O}$ is $((c_{ijk})_{i,j,k=1,\ldots,n})$;
$I=(\alpha,\beta)$, where $0\neq\alpha=\sum\limits_{i=1}^{n}a_i\omega_i$,
$\beta=\sum\limits_{i=1}^{n}b_i\omega_i$; $h_I$, a multiple-norm of $I$.
\ENSURE \quad TRUE iff $I$ is prime.

\STATE
Compute the factor ring $R=\mathcal{O}/I=\mathrm{MainAlgorithm}(\mathcal{O},I,h_I)$.
\STATE
Return $\mathrm{Is-Field}(R)$.
\end{algorithmic}
\end{algorithm}

The computational complexity of the field testing (Algorithm \ref{isfield}) is important
to deduce the time complexity of Algorithm $\ref{primeideal}$.

\begin{theorem}\label{sol1}\quad
Algorithm \rm{\ref{primeideal}} is correct and performs in
$O(\mathrm{M(log}^{15/2}h_I)+n^3\mathrm{B}(\mathrm{log}nT)+n^4\mathrm{B(log}h_I)\mathrm{log}h_I)$
bit operations, if all integers $c_{ijk},\;a_i,\;b_i,\;1\leq
i,j,k\leq n$, are bounded in absolute value by $T$.
\end{theorem}

\begin{proof}\quad Applying the relevant facts of \cite{AM} to Dedekind domains,
one may easily deduce that $I$ is a prime ideal if and only if $R=\mathcal{O}/I$ is a field.
Then the correctness follows immediately. Since we have
$|R|=d_1\cdot\ldots\cdot d_m=\mathcal{N}(I)|h_I$,
the time complexity of Algorithm \ref{primeideal} follows easily from
Theorem \ref{core1} and Lemma \ref{field}.
\end{proof}

\subsection{Deciding whether a given ideal is prime power}

The other application of the main algorithm is to decide whether
a given nonzero ideal is prime power in Dedekind domain of finite rank.
The local ring test Is-Local, which is stated in \cite{sta},
is crucial to deduce the following prime ideal power test.

\begin{algorithm}[htb]
\caption{Is-Prime-Power}\label{primepower}
\begin{algorithmic}[1]
\REQUIRE \quad $\mathcal{O}=\mathbb{Z}\omega_1\oplus\ldots\oplus\mathbb{Z}\omega_n$,
a multiplication table of $\mathcal{W}$ related to $\mathcal{O}$ is $((c_{ijk})_{i,j,k=1,\ldots,n})$;
$I=(\alpha,\beta)$, where $0\neq\alpha=\sum\limits_{i=1}^{n}a_i\omega_i$,
$\beta=\sum\limits_{i=1}^{n}b_i\omega_i$; $h_I$, a multiple-norm of $I$.
\ENSURE \quad TRUE iff $I$ is prime power.
\STATE
Compute the factor ring $R=\mathcal{O}/I=\mathrm{MainAlgorithm}(\mathcal{O},I,h_I)$.
\STATE
Return $\mathrm{Is-Local}(R)$.
\end{algorithmic}
\end{algorithm}

\begin{theorem}\label{sol2}\quad
Algorithm \rm{\ref{primepower}} is correct, and the time complexity is
$O(\mathrm{M(log}^{15/2}h_I)+n^3\mathrm{B}(\mathrm{log}nT)+n^4\mathrm{B(log}h_I)\mathrm{log}h_I)$
bit operations, where $T$ is as Theorem \ref{sol1}.
\end{theorem}

\begin{proof} Applying the relevant facts of \cite{AM} to Dedekind domains,
one may deduce that $I$ is a prime power if and only if
$R=\mathcal{O}/I$ is a local ring. Hence the correctness follows
immediately. Note that in Step 2 of Algorithm $\ref{primepower}$, the local ring test Is-Local(R) is deterministic and its computational
complexity is $O(\mathrm{M(log}^{15/2}p)+\mathrm{log}^4|R|)$
bit operations by the results of \cite{sta}, where $p=\mathrm{min}\{d_1,\ldots,d_m\}$.
According to Theorem \ref{core1}, the total time complexity of Algorithm $\ref{primepower}$ is obtained as asserted.
\end{proof}

\begin{remark}\quad
The input size of Algorithms \ref{primeideal} and \ref{primepower}
is $O(n^3\mathrm{log}T)$ bits. Similarly as Remark \ref{size}, one can verify that both Algorithm
\ref{primeideal} and Algorithm \ref{primepower} are polynomial time in the
input size. That is to say, our prime ideal test and prime ideal power test are
deterministic and polynomial-time.
\end{remark}

\section{Examples and implementation aspects of our algorithms}\label{Secimp}

In this section, we illustrate three different examples to show the implementation aspects of the main algorithm,
the prime ideal test and the prime ideal power test. The main algorithm (Algorithm \ref{main}) uses the
 algorithms for computing Hermite and Smith normal forms in \cite{hmc}, especially when computing the transforming matrices.
However, these algorithms are presented only with theoretical analysis of complexity, without empirical complexity and examples in \cite{hmc}.
The algorithms for fields testing and local rings testing are also stated with theoretical complexity
 in \cite{ADM} and \cite{sta} separately, without any implementation.
Hence, it is difficult to illustrate a large amount of empirical results in this paper,
even with the help of computer programming.
But we still calculate some concrete examples to illustrate each step of our algorithms explicitly.

\begin{example}\label{exm1}\quad
In: $\mathcal{O}=\mathbb{Z}[\delta]$, $I=(5,2+\delta)$, where $\delta=\sqrt[3]{2}$.
Out: $I$ is prime or not.

According to the prime ideal test (Algorithm \ref{primeideal}), we begin with computing
a basis representation of the factor ring $R=\mathcal{O}/I$. We implement each step of the
main algorithm (Algorithm \ref{main}) in the following way.

{\bf Step 1}\quad This step is the same as the one of Example \ref{exm}.

{\bf Step 2}\quad After few computations we may obtain the Hermite normal forms of $A^T$ and $B^T$ as follows,
which are denoted by $H_A$ and $H_B$ respectively,
$$H_A=
\left(
\begin{array}{ccc}
5 & 0  & 0  \\
0  & 5  & 0 \\
0  & 0 &  5
\end{array}\right);\quad
H_B=
\left(
\begin{array}{ccc}
10 & 2  &  6 \\
0  &  1 & 0 \\
0 & 0 & 1
\end{array}\right).$$

{\bf Step 3}\quad Similarly, we compute the Hermite normal form of $M=(H_A,H_B)$,
which is written as $(0\;H_M)$, and
$$H_M=
\left(
\begin{array}{ccc}
5 & 2  & 1  \\
0  & 1  & 0 \\
0  & 0 &  1
\end{array}\right).
$$

{\bf Step 4}\quad We continue to compute the Smith normal form $S$ of $H_M$ and
the related transforming matrices $U$, $V$, which are
$$V=
\left(
\begin{array}{ccc}
1 & -2  & -1  \\
0  & 1  & 0 \\
0  & 0 &  1
\end{array}\right);\quad
S=
\left(
\begin{array}{ccc}
5 & 0  & 0  \\
0  & 5  & 0 \\
0  & 0 &  5
\end{array}\right);$$
and $U=I_{3}$ is the identity matrix, such that $S=VH_MU$. Hence we get $d_1=5, d_2=d_3=1$.

{\bf Step 5}\quad Since $\mathrm{det}(S)=5>1$, we have $I\neq\mathcal{O}$.

{\bf Step 6}\quad It is easy to obtain the inverse of $V$ in this case:
$$V^{-1}=
\left(
\begin{array}{ccc}
1 & 2  & 1  \\
0  & 1  & 0 \\
0  & 0 &  1
\end{array}\right),$$
which leads to the same output as computing the matrix $\tilde{V}$ by Theorem \ref{core}. Hence we do not need to determine a multiple-norm $h_I$ in advance for this example. Thus we calculate all the integer matrices $A_k$ ($1\leq k\leq 3$) in the following:
$$A_1=
\left(
\begin{array}{ccc}
1 & 2  & 1  \\
2  & 4 & 4 \\
1  & 4 &  1
\end{array}\right);\quad
A_2=
\left(
\begin{array}{ccc}
0 & 1  & 0  \\
1  & 4 & 1 \\
0  & 1 &  2
\end{array}\right);\quad
A_3=
\left(
\begin{array}{ccc}
0 & 0  & 1  \\
0  & 1 & 2 \\
1  & 2 &  2
\end{array}\right).$$

{\bf Step 7}\quad After some computations we get
$$(t_{ij1})=
\left(
\begin{array}{ccc}
1 & 0  & 0 \\
0  & -5 & 0 \\
0  & 0 &  -5
\end{array}\right);\quad
(t_{ij2})=
\left(
\begin{array}{ccc}
0 & 1  & 0  \\
1  & 4 & 1 \\
0  & 1 &  2
\end{array}\right);\quad
(t_{ij3})=
\left(
\begin{array}{ccc}
0 & 0  & 1 \\
0  & 1 & 2 \\
1  & 2 &  2
\end{array}\right).$$

{\bf Step 8}\quad Since we only have $d_1=5>1$ and $l_{111}=\pi_1(t_{111})=1\pmod5=1$,
we obtain a basis representation of $R$, which is $R=\mathcal{O}/I=(1;5;(1))$.

Next we implement Step 2 of Algorithm \ref{primeideal}.
It is not difficult to deduce that Is-Field(R=(1;5;(1))) returns TRUE by Algorithm \ref{isfield}.
Hence $I$ is prime in this example.
\end{example}

\begin{example}\label{exmm}\quad
In: $\mathcal{O}=\mathbb{Z}[\theta]$, $I=(\theta-2)$, where $\theta=(1+\sqrt{-23})/2$.
Out: $I$ is prime power or not.

Note that, $\mathcal{O}=\mathbb{Z}\oplus\mathbb{Z}\theta$,
and the multiplication table of $\mathcal{W}=\{1,\theta\}$ is
$(c_{11k})_{k=1,2}=(1,0)$,
$(c_{12k})_{k=1,2}=(c_{21k})_{k=1,2}=(0,1)$,
$(c_{22k})_{k=1,2}=(-6,1)$.

According to the prime ideal power test (Algorithm \ref{primepower}), we start with computing
a basis representation of the factor ring $R=\mathcal{O}/I$. Each step of the
main algorithm (Algorithm \ref{main}) runs as follows.

{\bf Step 1}\quad By computation, the $2\times2$ integral matrix $A$ related to $I$ is
$$A=
\left(
\begin{array}{cc}
-2 & 1  \\
-6 & -1
\end{array}\right).$$
Since $I=(\theta-2)$ is a principal ideal, the matrix $B$ of the main algorithm
does not exist in this example.
Thus we only need to calculate the Hermite normal form of $A^T$ in Step 2 and Step 3 together,
which is
$$H_M=H_A=
\left(
\begin{array}{cc}
8 & 6  \\
0  & 1
\end{array}\right).$$

{\bf Step 4}\quad Now we compute the Smith normal form $S$ of $H_M$ and
the transforming matrices $U$ and $V$, which are
$$V=
\left(
\begin{array}{cc}
1 & -6  \\
0  & 1
\end{array}\right);\quad
S=
\left(
\begin{array}{cc}
8 & 0  \\
0  & 1
\end{array}\right)$$
and $U=I_{2}$ is the identity matrix, such that $S=VH_AU$. So we get $d_1=8, d_2=1$.

{\bf Step 5}\quad Since $\mathrm{det}(S)=8>1$, we have $I\subset\mathcal{O}$.

{\bf Step 6}\quad Similarly as Example \ref{exm1}, it is easy to obtain the inverse of $V$ at this time:
$$V^{-1}=
\left(
\begin{array}{cc}
1 & 6  \\
0  & 1
\end{array}\right).$$

Also, the integer matrices $A_k$ ($k=1, 2$) are computed by
$$A_1=
\left(
\begin{array}{cc}
1 & 6  \\
6 & 30
\end{array}\right);\quad
A_2=
\left(
\begin{array}{cc}
0 & 1  \\
1 & 13
\end{array}\right).$$

{\bf Step 7}\quad Under a few computations we may obtain
$$(t_{ij1})=
\left(
\begin{array}{cc}
1 & 0  \\
0  & -48
\end{array}\right);\quad
(t_{ij2})=
\left(
\begin{array}{cc}
0 & 1  \\
1  & 13
\end{array}\right).$$

{\bf Step 8}\quad Since only $d_1=8>1$ and $l_{111}=\pi_1(t_{111})=-35\pmod8=5$,
we obtain a basis representation of $R=\mathcal{O}/I$, which is $R=(1;8;(5))$.

Next we implement Step 2 of Algorithm \ref{primepower}. Simply applying the local ring test,
we have Is-Local(R=(1;8;(5))) returns TRUE.
Hence $I$ is prime power in this example.
Moreover, one may easily verify that Is-Field(R=(1;8;(5))) returns FALSE, that is to say,
$I$ is not a prime ideal.
\end{example}

\begin{example}\quad
In: $\mathcal{O}=\mathbb{Z}[\gamma]$, $I=(23)$, where
$\gamma^3=\gamma+1$.
Out: $I$ is prime power or not.

Note that, $\mathcal{O}=\mathbb{Z}\oplus\mathbb{Z}\gamma\oplus\mathbb{Z}\gamma^2$,
and the multiplication table of $\mathcal{W}=\{1,\gamma,\gamma^2\}$ is
$(c_{11k})_{k=1,2,3}=(1,0,0)$,
$(c_{12k})_{k=1,2,3}=(c_{21k})_{k=1,2,3}=(0,1,0)$,
$(c_{22k})_{k=1,2,3}=(c_{13k})_{k=1,2,3}=(c_{31k})_{k=1,2,3}=(0,0,1)$,
$(c_{23k})_{k=1,2,3}=(c_{32k})_{k=1,2,3}=(1,1,0)$,
$(c_{33k})_{k=1,2,3}=(0,1,1)$.

By the prime ideal power test (Algorithm \ref{primepower}), we may compute
a basis representation of the factor ring $R=\mathcal{O}/I$ firstly. Each step of the
main algorithm (Algorithm \ref{main}) is implemented as follows.

{\bf Step 1}\quad By computation, the $3\times3$ integral matrix $A$ related to $I$ is
$$A=
\left(
\begin{array}{ccc}
23 & 0 & 0  \\
0 & 23 & 0 \\
0 & 0 & 23
\end{array}\right)$$

The same as Example \ref{exmm}, the matrix $B$ does not appear in the main algorithm.
Since $A$ is already in Smith normal form and $I$ is principal,
we reach to Step 4 directly, that is,
the transforming matrices $U=V=I_{3}$ are the identity matrix and
$S=diag(23,23,23)$. Thus we have $d_1=d_2=d_3=23$.

{\bf Step 5}\quad Since $\mathrm{det}(S)=23^3=12167>1$, we get $I\subset\mathcal{O}$.

{\bf Step 6}\quad Because $V=I_{3}$ is the identity matrix here, it is easy to compute the integer matrices $A_k$ ($k=1, 2, 3$) in the following:
$$A_1=(c_{ij1})_{1\leq i, j\leq 3}=
\left(
\begin{array}{ccc}
1 & 0  & 0  \\
0  & 0 & 1 \\
0  & 1 &  0
\end{array}\right),$$

$$A_2=(c_{ij2})_{1\leq i, j\leq 3}=
\left(
\begin{array}{ccc}
0 & 1  & 0  \\
1  & 0 & 1 \\
0  & 1 &  1
\end{array}\right),$$

$$A_3=(c_{ij3})_{1\leq i, j\leq 3}=
\left(
\begin{array}{ccc}
0 & 0  & 1  \\
0  & 1 & 0 \\
1  & 0 &  1
\end{array}\right).$$

{\bf Step 7}\quad Under a few simple computations we may obtain
$$(t_{ij1})=
\left(
\begin{array}{ccc}
1 & 0  & 0  \\
0  & 0 & 1 \\
0  & 1 &  0
\end{array}\right);\quad
(t_{ij2})=
\left(
\begin{array}{ccc}
0 & 1  & 0  \\
1  & 0 & 1 \\
0  & 1 &  1
\end{array}\right);\quad
(t_{ij3})=
\left(
\begin{array}{ccc}
0 & 0  & 1  \\
0  & 1 & 0 \\
1  & 0 &  1
\end{array}\right).$$

{\bf Step 8}\quad Since $d_1=d_2=d_3=23>1$ and $l_{ijk}=\pi_k(t_{ijk})=t_{ijk}$ for all $i,j,k$.
We obtain the following basis representation of $R=\mathcal{O}/I$:
 $$R=(3;23,23,23;(1,0,0;0,0,1;0,1,0;0,1,0;1,0,1;0,1,1;0,0,1;0,1,0;1,0,1)).$$

Next we implement Step 2 of Algorithm \ref{primepower}. One may deduce that Is-Local(R=(3;23,23,23;$(t_{ijk})$)) returns FALSE by the local ring test in \cite{sta},
Hence $I$ is not prime power in this case.
\end{example}

One may verify that Is-Field(R=(3;23,23,23;$(t_{ijk})$)) also returns FALSE in the above example.

\section{Conclusion}\label{Seccon}
We have presented a deterministic polynomial-time algorithm for
computing a basis representation of the factor ring related to
a given ideal in Dedekind domain with finite rank. In addition,
we have also described the tests for deciding whether a nonzero ideal is prime
and whether it is prime power as the important
applications of the main algorithm. The prime ideal test and the prime ideal
power test are proven to be deterministic and polynomial-time complexity.
All the algorithms use $h_I$ which is a multiple-norm of $I$, a
smaller $h_I$ makes these algorithms more efficient.

If an integral basis of the ring of algebraic integers of a number
field is known, then a deterministic polynomial-time
algorithm for testing the primality of ideals in this ring can be
deduced from our prime ideal test. It is natural to
ask whether there exists a deterministic polynomial-time algorithm
for computing the factor rings or testing the primality of ideals in a general Dedekind domain,
not necessarily of finite rank. At this time, the corresponding
factor ring need not be a finite ring, such as
$\mathcal{O}=\mathbb{Q}[X],\;I=(X)$ and
$\mathcal{O}/I\cong\mathbb{Q}$ is infinite. Hence the current method
based on the main algorithm will not work any more. We are looking
forward to finding a new method in the future work.

\Acknowledgements{The authors would like to thank Dr. Yunling Kang for his helpful discussion.
This work was supported by the NSFC (Grant No. 11601202), the NNSF of China (Grant No. 11471314),
the NSFC (Grant No. 11401312), the NSF of the Jiangsu Higher Education
Institutions (Grant No. 14KJB110012), the high-level talent scientific research foundation of Jinling
Institute of Technology grant jit-b-201527, and the National
Center for Mathematics and Interdisciplinary Sciences, CAS.}


\end{document}